\documentclass[12pt,reqno]{amsart}
\usepackage{fullpage}
\usepackage{times}
\usepackage[utf8]{inputenc}
\usepackage[english]{babel}
\usepackage[mathscr]{euscript}
\usepackage{amsmath,amssymb,amsthm, graphicx, enumitem, xcolor}
\usepackage[colorlinks=true,citecolor=blue,linkcolor=blue,urlcolor=blue]{hyperref}
\newtheorem{theorem}{Theorem}[section]
\newtheorem{proposition}[theorem]{Proposition}
\newtheorem{lemma}[theorem]{Lemma}

\newtheorem{question}[theorem]{Question}

\theoremstyle{definition}

\newcommand{\R}{\mathbb{R}} 

\newcommand{\Z}{\mathbb{Z}}
\newcommand{\T}{\mathbb T}

\title{A resolution of the Erd\H{o}s similarity problem in the large}
\author{Chi Hoi Yip}
\address{Department of Mathematics, Hong Kong University of Science and Technology, Clear Water Bay, Hong Kong}
\email{machyip@ust.hk}
\begin{document}
\maketitle

\begin{abstract}
We show that for every unbounded set $A\subset\mathbb R$ and every $0<\varepsilon<1$, there exists a closed set $S\subset\mathbb R$ that has measure at least $1-\varepsilon$ in every unit interval and contains no nontrivial affine copy of $A$. Moreover, $S$ can be chosen so that $|S\cap[x,x+1]|\to1$ as $|x|\to\infty$. This resolves the Erd\H{o}s similarity problem in the large. 
\end{abstract}

\section{Introduction}
A standard application of the Lebesgue density theorem shows that for every finite set $A \subset \mathbb{R}$ and a set $S \subset \mathbb{R}$ with positive (Lebesgue) measure, $S$ contains a nontrivial affine copy of $A$; that is, there exist $\lambda \in \R\setminus \{0\}$ and $t \in \R$ such that \[\lambda A+t=\{\lambda a +t: a\in A\} \subseteq S.\] This was first proved by Steinhaus \cite{S20}. In 1974, motivated by this statement, Erd\H{o}s \cite{Erdos74} conjectured that all infinite sets fail to satisfy this property, that is, for each infinite set $A \subset \mathbb{R}$, there is a set $S \subset \mathbb{R}$ with positive measure such that $S$ contains no nontrivial affine copy of $A$. This is known as the \emph{Erd\H{o}s similarity problem} or the \emph{Erd\H{o}s similarity conjecture}. 

As Bourgain \cite{Bourgain87} pointed out, Erd\H{o}s' question was presumably motivated by the search for an infinite analogue of the finite-pattern phenomenon guaranteed by Szemer\'edi's theorem~\cite{Szemeredi75} and its extension by Furstenberg and Katznelson~\cite{FK78}. Indeed, Szemer\'edi's theorem asserts that a set of integers with positive upper density contains arbitrarily long finite arithmetic progressions, while Steinhaus' theorem gives the continuous analogue that every measurable set of positive Lebesgue measure contains affine copies of every finite configuration.  The Erd\H{o}s similarity conjecture asks whether this phenomenon must fail for infinite configurations.

It is known~\cite{Sve01} that to prove the Erd\H{o}s similarity conjecture, it suffices to consider the case when $A$ is an infinite sequence decreasing to $0$. The conjecture has been proved for several families of sequences; see for example Falconer~\cite{Falconer84}, Eigen \cite{Eigen85}, Bourgain~\cite{Bourgain87},  Kolountzakis~\cite{Kolountzakis97}, Humke and Laczkovich~\cite{HL98}. We refer the reader to \cite{Sve01, JLM} for comprehensive surveys of the conjecture. 

Recent developments on the Erd\H{o}s similarity problem include results concerning bi-Lipschitz embeddings \cite{FLX24}, Cantor sets \cite{Kolountzakis23,SY26}, and sumsets of infinite sequences \cite{MICKRY26}. Very recently, OpenAI~\cite{OpenAI26} announced a proof of the conjecture for geometric progressions. See Iosevich and Yavicoli~\cite{IY26} for an exposition of the proof, and Jiang and Xi~\cite{JX26} for an extension to null sequences whose successive ratios are bounded away from zero. Nevertheless, the Erd\H{o}s similarity conjecture remains open in general.

In this paper, we study the Erd\H{o}s similarity problem in the large. For a Lebesgue measurable set $E$, we use $|E|$ to denote its Lebesgue measure. Following~\cite{GMY, IKMRY}, for $0<\varepsilon<1$, we say that a Lebesgue measurable set $S\subset\R$ is \emph{$(1-\varepsilon)$-large} if
\begin{equation}\label{eq:large}
 |S\cap I|\ge1-\varepsilon
 \qquad\text{for every unit interval $I\subset\R$.}
\end{equation}
Motivated by \cite{BKM}, Kolountzakis and
Papageorgiou~\cite[Question~1]{KP} (see also \cite[Question 1.2]{GMY}) asked the following question:
\begin{question}\label{q}
Does there exist an increasing sequence $A=\{a_n:n\ge1\}$ with $a_n\to\infty$ and some $0<\varepsilon<1$ such that every
$(1-\varepsilon)$-large set $S\subset\R$ contains a nontrivial affine copy of $A$?
\end{question}

This question is known as the \emph{Erd\H{o}s similarity problem in the large}. Question~\ref{q} can also be viewed as a continuous, infinite counterpart of Szemer\'edi's theorem,
asking whether ``sufficiently large'' sets must contain affine copies of some infinite configuration.

Question~\ref{q} has been extensively studied in several recent works; see \cite{BKM, BGKMW,GMY,IKMRY,JLM,KP}. For each $\varepsilon>0$, Bradford, Kohut, and Mooroogen~\cite{BKM} explicitly constructed $(1-\varepsilon)$-large sets containing no infinite arithmetic progression.  Using the probabilistic method, Kolountzakis and Papageorgiou~\cite{KP} answered Question~\ref{q} in the negative for separated sequences of subexponential growth. Using tools from metric number theory, Gao, Mooroogen, and Yip~\cite{GMY} answered the question in the negative for certain sequences of arbitrarily rapid growth. More recently, Iosevich, Kulkarni, Mora Cu\'ellar,
Rojas Aravena, and Yavicoli~\cite{IKMRY} constructed
periodic $(1-\varepsilon)$-large sets avoiding affine
copies of sets supporting Rajchman measures (that is, measures whose Fourier transform vanishes at infinity). However, these results did not address arbitrary unbounded sets.

In this paper, we answer Question~\ref{q} negatively for all unbounded subsets of $\R$.
\begin{theorem}\label{thm:main}
Let $A\subset\R$ be an unbounded set, $0<\varepsilon<1$, and $R_0>0$.
There exists a closed, symmetric, $(1-\varepsilon)$-large set
$S\subset\R$ containing $[-R_0,R_0]$ such that
\begin{equation}\label{eq:limit-density}
 \lim_{|x|\to\infty}|S\cap[x,x+1]|=1,
\end{equation}
and
\begin{equation}\label{eq:infinite-omissions}
 (\lambda A+t)\setminus S \text{ is an infinite set}
\end{equation}
for all $\lambda\ne0$ and $t\in\R$.
\end{theorem}
We remark that the rate of convergence in condition~\eqref{eq:limit-density} cannot, in general, be prescribed independently of \(A\). See Section~\ref{sec3}.

\section{Proof of the main theorem}
Let $\T=\R/\Z$. For $u\in \R$, let $\|u\|$ be the distance from $u$ to the nearest integer.

We first establish a local version of Theorem~\ref{thm:main}.

\begin{lemma}\label{lem:finite}
Let $A\subset(0,\infty)$ be unbounded. Given $0<\delta<1$, $0<c\leq C$, and $T,R>0$, there exist a finite set $F\subset A$
and a bounded open set $U\subset(R,\infty)$ such that
$\R\setminus U$ is $(1-\delta)$-large and
\begin{equation}\label{eq:finite-hit}
 (\lambda F+t)\cap U\ne\varnothing
 \qquad(c\le\lambda\le C,\ |t|\le T).
\end{equation}
\end{lemma}

We first give a simple probabilistic proof of Lemma~\ref{lem:finite} to illustrate the main idea. After proving the main theorem, we show how to make the construction deterministic using quantitative
equidistribution.

\subsection{A probabilistic proof of Lemma~\ref{lem:finite}}

For $a\in A$, set
\[
I(a)=(ca-T-1,\,Ca+T+1).
\]
Note that for every $c\le\lambda\le C$ and $|t|\le T$, the point
$\lambda a+t$ lies in the interval $I(a)$.

Set $h=\delta/16$ and choose an integer $J\ge1/h$. Choose $a_0\in A$ with $I(a_0)\subset(R,\infty)$,
and then choose an integer $k$ so large that
\begin{equation}\label{eq:J}
 J|I(a_0)|(1-2h)^k<\frac{\delta}{4}.
\end{equation}
Since $A$ is unbounded, we can choose
$a_0<a_1<\cdots<a_k$ from $A$ so that the intervals
$I(a_i)$, $0\le i\le k$, have pairwise distances greater
than $2$. 

For $1\le i\le k$, define
\[
 r_i=\frac{a_i}{a_i-a_0}>1.
\]
Choose $L>0$ sufficiently large that
\begin{equation}\label{eq:L}
 L(r_i-1)>\frac{8}{\delta}\qquad(1\le i\le k).
\end{equation}

Let $\theta_1,\ldots,\theta_k$ be independent random variables, each uniformly distributed on $\T$. For fixed $y\in I(a_0)$ and $0\le j<J$, the events
\[
\left\|Lr_i y-\theta_i-\frac{j}{J}\right\|>h,\qquad 1\le i\le k,
\]
are independent and each has probability $1-2h$.
Consequently, the set
\begin{equation}\label{eq:U0}
 U_0=\bigcup_{j=0}^{J-1}
 \left\{y\in I(a_0):
 \left\|Lr_i y-\theta_i-\frac{j}{J}\right\|>h
 \text{ for all }1\le i\le k\right\}
\end{equation}
satisfies, by the union bound, Fubini's theorem, and estimate~\eqref{eq:J}, 
\[
 \mathbb E|U_0|
 \le\sum_{j=0}^{J-1}\int_{I(a_0)}(1-2h)^k\,dy
 =J|I(a_0)|(1-2h)^k<\delta/4.
\]
Thus, we can find $\theta_1,\ldots,\theta_k\in \T$ for which
\begin{equation}\label{eq:center-small}
 |U_0|<\delta/4.
\end{equation}

For $1\le i\le k$, define
\begin{equation}\label{eq:Ui}
 U_i=\left\{y\in I(a_i):
 \|L(r_i-1)y-\theta_i\|<2h\right\}.
\end{equation}
Observe that for every $x\in\R$, by the choice of $h$ and estimate~\eqref{eq:L}, we have
\begin{equation}\label{eq:measure}
 |U_i\cap[x,x+1]|
 \le4h+\frac{2}{L(r_i-1)}<\frac{\delta}{2}.
\end{equation}
Set
\[
 F=\{a_0,\ldots,a_k\},\qquad U=\bigcup_{i=0}^k U_i.
\]
The set $U$ is bounded and open. Since a unit interval can
intersect at most one of the intervals $I(a_i)$, the bounds
\eqref{eq:center-small} and \eqref{eq:measure} show that
$\R\setminus U$ is $(1-\delta)$-large.

It remains to verify that $U$ intersects every affine copy of $F$ under consideration. Fix $c\le\lambda\le C$ and $|t|\le T$, and set $y_i=\lambda a_i+t\in I(a_i)$ for $0\leq i\leq k$. If $y_0\in U_0$, then we are done. Next assume $y_0\notin U_0$. Choose $0\leq j\leq J-1$ so that \[\left\|Lt-\frac{j}{J}\right\|\leq \frac{1}{2J}\leq \frac{h}{2}.\] By the definition of $U_0$ in equation~\eqref{eq:U0}, there is $1\le i\le k$ with
\begin{equation}\label{eq:center-near}
 \left\|Lr_i y_0-\theta_i-\frac{j}{J}\right\|\le h.
\end{equation}
Since $y_0=\lambda a_0+t$ and $y_i=\lambda a_i+t$, we have
\begin{equation}\label{eq:key-identity}
 r_i y_0-(r_i-1)y_i=\frac{a_i}{a_i-a_0}(\lambda a_0+t)
 -\frac{a_0}{a_i-a_0}(\lambda a_i+t)=t.
\end{equation}
Using estimate~\eqref{eq:center-near} and equation~\eqref{eq:key-identity}, we get
\[
 \begin{aligned}
 \|L(r_i-1)y_i-\theta_i\|
 &=\|Lr_i y_0-\theta_i-Lt\|\\
 &\le \left\|Lr_i y_0-\theta_i-\frac{j}{J}\right\|+\left\|Lt-\frac{j}{J}\right\|\\
 &<2h.
 \end{aligned}
\]
Thus $y_i\in U_i$ by the definition of $U_i$ in equation~\eqref{eq:Ui}, completing the proof.

\subsection{Proof of Theorem~\ref{thm:main}}
By replacing $A$ with an unbounded subset of $A$ or $-A$,
we may assume that $A\subset(0,\infty)$ is unbounded.

We construct bounded open sets $U_s$ inductively.
Set $R_1=R_0+3$. Suppose that $R_s$ has been chosen.
Applying Lemma~\ref{lem:finite} with
\[
c=\frac{1}{s},\qquad C=s,\qquad T=s,\qquad
R=R_s,\qquad \delta=\delta_s=\frac{\varepsilon}{2^{s+1}},
\]
we obtain a bounded open set $U_s\subset(R_s,\infty)$
such that
\[
|U_s\cap[x,x+1]|\le\delta_s
\qquad(x\in\R),
\]
and
\[
(\lambda A+t)\cap U_s\ne\varnothing
\qquad(s^{-1}\le\lambda\le s,\ |t|\le s).
\]
We then set
\[
R_{s+1}=\sup \{x: x\in U_s\}+3
\]
and continue the construction.

Now define
\[
U=\bigcup_{s=1}^{\infty} \bigl(U_s\cup(-U_s)\bigr),
\qquad S=\R\setminus U.
\]
Then $S$ is closed, symmetric, and contains
$[-R_0,R_0]$. By construction, the sets $U_s$
and their reflections are pairwise separated by
a distance greater than $2$. Thus every unit
interval intersects at most one of them, and
\[
|U\cap[x,x+1]|
\le\sup_{s\ge1}\delta_s<\varepsilon
\qquad(x\in\R).
\]
Hence $S$ is $(1-\varepsilon)$-large. Moreover, since $R_s \to \infty$ and $\delta_s \to 0$ as $s\to \infty$, we have
\[
\lim_{|x|\to\infty}|S\cap[x,x+1]|=1.
\]

Finally, we verify condition~\eqref{eq:infinite-omissions}. Fix $\lambda\ne0$ and $t\in\R$. By the symmetry of $U$, we may assume that $\lambda>0$. For all sufficiently large $s$, we have
$s^{-1}\le\lambda\le s$ and $|t|\le s$. Thus, by construction,
\[
(\lambda A+t)\cap U_s\ne\varnothing.
\]
Since the sets $U_s$ are pairwise disjoint,
$(\lambda A+t)\setminus S$ contains infinitely
many points, as required. This finishes the proof of the theorem.

\subsection{Making the proof deterministic}

The only randomness in the above proof of  Theorem~\ref{thm:main} comes from choosing the independent shifts $\theta_1,\ldots,\theta_k$ in the proof of Lemma~\ref{lem:finite}. We explain
why they can all be taken to be zero if the $a_i$ are chosen sufficiently
sparsely. The argument uses a quantitative form of the classical
Kronecker--Weyl equidistribution theorem; see also
\cite{GM16,Mak25} for related quantitative approximation results.

We record a special case of the continuous Erd\H{o}s--Tur\'an--Koksma inequality below, following \cite[Theorem~2.77]{DT97}. For $\mathbf{m}=(m_1,m_2,\ldots, m_k)\in \Z^k$, we write
\[
\|\mathbf{m}\|_\infty=\max \{|m_i|: 1\leq i \leq k\}, \qquad 
w(\mathbf{m})=\prod_{i=1}^k\max\{1,|m_i|\}.
\]

\begin{lemma}
\label{lem:continuous-ETK}
Let $I\subset\R$ be a bounded interval of positive length and let
$f=(f_1,\ldots,f_k):I\to\R^k$ be continuous. For every axis-parallel
box $B\subset[0,1)^k$ and every positive integer $H$,
\begin{equation}\label{eq:continuous-ETK}
\begin{aligned}
&\left|\frac1{|I|}\int_I\mathbf{1}_B(\{f(y)\})\,dy-|B|\right|\\
&\qquad\le\left(\frac32\right)^k
\left(\frac2{H+1}
+\sum_{\substack{\mathbf{m}\in\Z^k\\0<\|\mathbf{m}\|_\infty\le H}}
\frac1{w(\mathbf{m})}
\left|\frac1{|I|}\int_I
 e\bigl(\mathbf{m}\cdot f(y)\bigr)\,dy\right|\right),
\end{aligned}
\end{equation}
Here $e(t)=e^{2\pi it}$ and fractional parts are taken
coordinatewise. 
\end{lemma}

We apply Lemma~\ref{lem:continuous-ETK} with
$f_i(y)=Lr_i y-u$ and $B=(h,1-h)^k$, where $L>0$ and $h \in (0,1/4)$. Set
\[
\Delta_H=
\min_{\substack{\mathbf{m}\in\Z^k\setminus\{0\}\\|m_i|\le H}}
\left|\sum_{i=1}^k m_i r_i\right|.
\]
If $\Delta_H>0$, direct integration gives, for every nonzero $\mathbf{m}$
with $\|\mathbf{m}\|_\infty\le H$,
\[
\left|\int_I e\bigl(\mathbf{m}\cdot f(y)\bigr)\,dy\right|
\le\frac1{\pi L\Delta_H}.
\]
Since
\[
\sum_{\substack{\mathbf{m}\in\Z^k\\0<\|\mathbf{m}\|_\infty\le H}}
\frac1{w(\mathbf{m})}
\le\left(1+2\sum_{j=1}^H\frac1j\right)^k
\le(3+2\log H)^k,
\]
we obtain the estimate
\begin{equation}\label{eq:det-discrepancy}
\begin{aligned}
&\bigg|\big|\{y\in I:
\|Lr_i y-u\|>h\ \text{for all }1\le i\le k
\}\big|-|I|(1-2h)^k\bigg|\\
&\qquad\le\left(\frac32\right)^k
\left(\frac{2|I|}{H+1}
+\frac{(3+2\log H)^k}{\pi L\Delta_H}\right),
\end{aligned}
\end{equation}
uniformly in $u\in\R$.

We now return to Lemma~\ref{lem:finite}, retaining $h$, $J$, and $a_0$
from its proof. Choose $k$ so large that
\begin{equation}\label{eq:J1}
J|I(a_0)|(1-2h)^k<\frac{\delta}{8},
\end{equation}
and then $H$ so large that
\begin{equation}\label{eq:J2}
\left(\frac32\right)^k\frac{2J|I(a_0)|}{H+1}
<\frac\delta{16}.
\end{equation}
Since $A$ is unbounded, we can choose $a_1,\ldots,a_k \in A$ successively,
still keeping the intervals $I(a_i)$ separated, so that, with
$b_i=r_i-1$ and $\kappa=1/(4(H+1))$,
\[
0<b_1<\kappa,\qquad
0<b_{i+1}<\kappa b_i\quad(1\le i<k).
\]

These inequalities imply $\Delta_H\ge\tfrac34 b_k>0$.
Indeed, let $\mathbf{m}=(m_1,\ldots,m_k)\in\mathbb Z^k\setminus\{0\}$
with $\|\mathbf{m}\|_\infty\le H$. Since $r_i=1+b_i$, we have
\[
\sum_{i=1}^k m_i r_i
=\sum_{i=1}^k m_i+\sum_{i=1}^k m_i b_i.
\]
If $\sum_{i=1}^{k} m_i\ne0$, then $|\sum_i m_i|\ge1$, while
\[
\left|\sum_{i=1}^k m_i b_i\right|
\le H\sum_{i=1}^k b_i
<\frac{H\kappa}{1-\kappa}<\frac14.
\]
Consequently, $|\sum_{i=1}^k m_i r_i|>\frac34\ge\tfrac34 b_k$. Otherwise, let $\ell$ be the smallest index with $m_\ell\ne0$. Since $\sum_{i=1}^k m_i=0$, we have 
\[
\begin{aligned}
\left|\sum_{i=1}^k m_i r_i\right|
=\left|\sum_{i=\ell}^k m_i b_i\right|
\ge b_\ell-H\sum_{i>\ell}b_i
>b_\ell-\frac{H\kappa}{1-\kappa}b_\ell
>\frac34 b_\ell
\ge\frac34 b_k.
\end{aligned}
\]
Thus $\Delta_H\ge\tfrac34 b_k$, as claimed.

Finally, choose $L$ large enough to satisfy inequality~\eqref{eq:L} and
\begin{equation}\label{eq:J3}
\left(\frac32\right)^k
\frac{J(3+2\log H)^k}{\pi L\Delta_H}<\frac\delta{16}.
\end{equation}
Taking $\theta_1=\cdots=\theta_k=0$ in equation~\eqref{eq:U0}, the union bound
and inequalities~\eqref{eq:det-discrepancy},~\eqref{eq:J1},~\eqref{eq:J2}, and~\eqref{eq:J3} together give
\[
|U_0|\leq J|I(a_0)|(1-2h)^k+J\left(\frac32\right)^k
\left(\frac{2|I(a_0)|}{H+1}
+\frac{(3+2\log H)^k}{\pi L\Delta_H}\right)<\frac\delta8+\frac\delta{16}+\frac\delta{16}
=\frac\delta4.
\]
Thus inequality~\eqref{eq:center-small} holds. The rest
of the proof of Lemma~\ref{lem:finite} is unchanged, yielding a
deterministic construction of $U$.

\section{A restriction on the rate of convergence}\label{sec3}
By Theorem~\ref{thm:main}, we can choose $S$ so that
\[
|S\cap[x,x+1]|\to1
\qquad\text{as }|x|\to\infty.
\]
However, the rate of convergence cannot, in general,
be prescribed independently of $A$. This follows from a simple
probabilistic argument.

\begin{proposition}\label{prop:rate}
Let $\eta:[1,\infty)\to(0,\infty)$ be a decreasing function such that
$\eta(x)\to0$ as $x\to\infty$. Suppose that
$1\le a_1<a_2<\cdots$ and
\[
 \sum_{n=1}^{\infty}\eta(a_n)<\infty.
\]
If a measurable set $S\subset\R$ satisfies
\[
 |[x,x+1]\setminus S|\le\eta(x)\qquad(x\ge1),
\]
then $S$ contains a positive dilate of $\{a_n:n\ge1\}$.
\end{proposition}

\begin{proof}
Put $U=\R\setminus S$, and let $\lambda$ be uniformly
distributed on $[B,2B]$, where $B\ge1$. For each $n$,
\[
 \mathbb P(\lambda a_n\in U)
 =\frac{|U\cap[Ba_n,2Ba_n]|}{Ba_n}.
\]
By the assumption on $S$ and the monotonicity of $\eta$,
\[
 \mathbb P(\lambda a_n\in U)
 \le\left(1+\frac1{Ba_n}\right)\eta(Ba_n)
 \le 2\eta(Ba_n).
\]
For each fixed $n$, this bound tends to $0$ as $B\to\infty$. Since $\sum_{n=1}^{\infty}\eta(a_n)<\infty$, dominated convergence gives
\[
 \sum_{n=1}^{\infty}\mathbb P(\lambda a_n\in U)\longrightarrow0.
\]
as $B\to \infty$. In particular, for sufficiently large $B$, the sum is less than $1$. Therefore, there is then a $\lambda\in[B,2B]$ for which $\lambda a_n\notin U$ for every $n$, as claimed.
\end{proof}

Indeed, given any decreasing function \(\eta(x)\to0\), we may choose an increasing sequence \(a_n\to\infty\) so rapidly that \(\eta(a_n)\le2^{-n}\). Proposition~\ref{prop:rate} then shows that no set avoiding all nontrivial affine copies of \(\{a_n:n\ge1\}\) can satisfy the prescribed rate of convergence.

A similar argument gives a restriction on the measure of the  complement of a set containing no nontrivial affine copy of $A$. Indeed, suppose that $A=\{a_n:n\ge1\}$ satisfies $a_n\ge1$ and
\[
\sum_{n=1}^{\infty}\frac1{a_n}<\infty.
\] 
If $U\subset\R$ is measurable with $|U|<\infty$, then for sufficiently large $B$, choosing $\lambda$ uniformly from $[B,2B]$ yields
\[
\sum_{n=1}^{\infty}
\mathbb P(\lambda a_n\in U)
\le\frac{|U|}{B}
\sum_{n=1}^{\infty}\frac1{a_n}<1.
\]
Thus $S=\R\setminus U$ contains a positive dilate of $A$. In particular, any set avoiding all nontrivial affine copies of $\{b^n:n\ge1\}$, where $b>1$, must have a complement of infinite measure.

\section*{Acknowledgments}
The author thanks Xiang Gao and Joshua Zahl for many valuable suggestions. The author also thanks Alex Burgin, Yuveshen Mooroogen, and Malabika Pramanik for helpful discussions at the early stage of the project. 

The author acknowledges the use of ChatGPT by OpenAI to assist with the exposition of the manuscript. The mathematical ideas are due to the author.

\end{document}